\documentclass[12pt,fleqn]{article}
\usepackage{amssymb,amsmath,eucal,amsthm}
\def\ds{\displaystyle}
\def\ov{\overline}
\def\bce\{\begin{center}
\def\ece{\end{center}}
\def\beq\{\begin{equation}
\def\eeq{\end{equation}}
\def\bsp\{\begin{split}
\def\esp{\end{split}}

\newtheorem{proposition}{Proposition}[section]
\theoremstyle{definition}

\theoremstyle{remark}

\renewcommand{\proofname}{\bf Proof}

\title{The geometry of fractional osculator bundle of higher order and applications}
\author{Ion Doru Albu$^a$, Mihaela Neam\c {t}u$^b$, Dumitru Opri\c {s}$^c$}
\date{ }
\begin{document}
\maketitle

\begin{tabular}{cccccccc}
\scriptsize{$^{a}$ Department of Mathematics, Faculty of
Mathematics and Informatics,
West University of Timi\c soara,}\\
\scriptsize{Bd. V. Parvan, nr. 4, 300223, Timi\c soara, Romania, e-mail: albud@math.uvt.ro,}\\
\scriptsize{$^{b}$Department of Economic Informatics and Statistics, Faculty of Economics,West University of Timi\c soara,}\\
\scriptsize{Str. Pestalozzi, nr. 16A, 300115, Timi\c soara, Romania, e-mail:mihaela.neamtu@fse.uvt.ro,}\\
\scriptsize{$^{c}$ Department of Applied
Mathematics, Faculty of Mathematics, West University of Timi\c soara,}\\
\scriptsize{Bd. V. Parvan, nr. 4, 300223, Timi\c soara, Romania, e-mail: opris@math.uvt.ro.}\\

\end{tabular}

\begin{abstract}Using the reviewed Riemann-Liouville fractional derivative we
define the bundle ${\mathop{E}\limits^{\alpha k}}=Osc^{\alpha
k}(M)$ and highlight geometrical structures with a geometrical
character. Also, we introduce the fractional osculator Lagrange
space of k order and the main structures on it. The results are
applied at the k order fractional prolongation of Lagrange,
Finsler and Riemann fractional structures.
\end{abstract}

\vspace{3mm}


\vspace{2mm}

\section{Introduction}
\hspace{0.5cm} It is known that the operators of integration and
derivation have geometrical and physical interpretations and they
were used in the modelation of problems from different domains.
The use of reviewed Liouville-Riemann integration and derivation
operators lead to fractional integration and derivation. The
geometrical and physical interpretation is suggested by the
Stieltjes integral and it was done by I. Podlubny [7]. There is a
vast bibliography which contains the properties of fractional
integral and derivative and the analysis of the processes which
are modeled with their help [2], [5], [8].

A lot of models which use the fractional derivative are defined on
an open set in ${{\mathbb R}}^n$. In this paper we present the
fractional derivative taking into account the geometrical
character, namely the behavior of associated objects under a
change of local chart.

The outline of this paper is as follows. In Section 2 we describe
the reviewed fractional derivative on ${\mathbb R}$ using [2],
[5], the fractional osculator bundle $T^\alpha(M)$ on a
differentiable manifold and the behavior of introduced objects
under a change of local chart. In Section 3, we define the
fractional osculator bundle of k order using the method presented
by R. Miron in [6]. We introduce: the Liouvile fractional vector
fields, the $\alpha k$-fractional spray and the fractional
nonlinear connection. We prove that these objects have a
geometrical character. Our findings are analogous with R. Miron's
results for the fractional case. In Section 4 we describe the
fractional Euler-Lagrange equations for fractional osculator
Lagrange spaces of superior order. The results are applied for the
k-order fractional bundle prolongation of Lagrange, Finsler and
Riemann structures.

The main results from the present paper were used in [3] and [4]
for the study of some fractional geometrical structures and they
will permit the study of other structures of this type.

\section{The fractional derivative on R. The fractional osculator bundle on the differentiable manifold.}

\subsection{The fractional derivative on R}

\hspace{0.5cm} Let $f:[a,b]\rightarrow \mathbb R$ be a derivable
function and $\alpha\in \mathbb R$, $\alpha>0$. The functions:
\begin{equation*}
\begin{split}
&{}_aD{}^\alpha_tf(t)=\frac{1}{\Gamma(m-\alpha)}(\frac{d}{dt})^m\int^t_a(t-s)^{m-\alpha-1}(f(s)-f(a))ds\\
&{}_tD{}^\alpha_bf(t)=\frac{1}{\Gamma(m-\alpha)}(-\frac{d}{dt})^m\int^b_t(t-s)^{m-\alpha-1}(f(s)-f(b))ds,
\end{split}
\end{equation*} are called the left respectively the right Liouville Riemann
fractional derivatives of function f, where $m\in \mathbb N^*$
with $m-1\leq\alpha<m$ and $\Gamma$ is Euler Gamma function.

In general, the operators ${}_aD{}^\alpha_t$, ${}_tD{}^\alpha_b$
do not satisfy semigroupal properties with respect to the
concatenation operation. Thus, we define the derivative operators
on the function spaces where the semigroupal properties hold.

The functions:
\begin{equation*} \begin{split}
&D{}^\alpha_tf(t)\!=\!\frac{1}{\Gamma(m\!-\!\alpha)}(\frac{d}{dt})^m\int^t_{-\infty}(t-s)^{m\!-\!\alpha\!-\!1}(x(s)\!-\!x(0))ds,\quad 0\in (\!-\!\infty,t)\\
&{}^*D{}^\alpha_tf(t)\!=\!\frac{1}{\Gamma(m\!-\!\alpha)}(\!-\frac{d}{dt})^m\int^{\infty}_t(s\!-\!t)^{m\!-\!\alpha\!-\!1}(x(s)\!-\!x(0))ds,
\quad 0\in (t,\infty)
\end{split} \end{equation*}
are called the left, respectively the right fractional derivative
of $\alpha$ order for function f.

If $\ov{supp}(f)=C(a,b),$ then $D_t^\alpha f(t)=_aD^\alpha_tf(t)$,
$^*D_t^\alpha f(t)=_tD^\alpha_bf(t)$.

We define the seminorms:

\begin{equation*}|x|_{J^\alpha_L(R)}=||D^\alpha_t||_{L^2(R)}, \quad
|x|_{J^\alpha_R(R)}=||^*D^\alpha_t||_{L^2(R)}\end{equation*} and
the norms:
\begin{equation*}||x||_{J^\alpha_L(R)}=(||x||^2_{L^2(R)}+|x|^2_{J^\alpha_L(R)}), \quad
||x||_{J^\alpha_R(R)}=(||x||^2_{L^2(R)}+|x|^2_{J^\alpha_L(R)}),\end{equation*}
where $J^\alpha_L(R)$, respectively $J^\alpha_R(R)$ denotes the
closure of $C_0^\infty(R)$ with respect to
$||\cdot||_{J^\alpha_L(R)}$, respectively
$||\cdot||_{J^\alpha_R(R)}$.

From the above definitions we have the following [2]:

\begin{proposition}\label{prop1}1. Let $I\subset R$ and
$J^\alpha_{L,0}(I)$, $J^\alpha_{R,0}(I)$ be the closure of
$C_0^\infty(I)$ in accordance with the respective norms. Then, for
any $f\in J^\beta_{L,0}(I)$, $0<\alpha<\beta$, respectively for
any $f\in J^\beta_{R,0}(I)$, $0<\alpha<\beta$, the relation
\begin{equation*}
D^\beta_tf(t)=D^\alpha_tD_t^{\beta-\alpha}f(t),
\end{equation*} respectively, the relation
\begin{equation*}
^*D^\beta_tf(t)=^*D^\alpha_t{}^*D_t^{\beta-\alpha}f(t)
\end{equation*}
holds;

2. If $\underset{n\rightarrow\infty}{lim}\alpha_n=p\in\mathbb
N^*,$ then
\begin{equation*}
\underset{n\rightarrow\infty}{lim}(D^{\alpha_n}_tf(t))=D^p_tf(t),\quad
\underset{n\rightarrow\infty}{lim}(^*D^{\alpha_n}_tf(t))=^*D^{p}f(t);
\end{equation*}

3. (i) If $f(t)=c$, $t\in [a,b]$, then $D^\alpha_tf(t)=0$;

(ii) If $f_1(t)=t^\gamma$, $t\in [a,b]$, then
$D^\alpha_tf_1(t)=\ds\frac{t^{\gamma-\alpha}\Gamma(1+\gamma)}{\Gamma(1+\gamma-\alpha)}$;

4. If $f_1$, $f_2$ are analytical functions on $[a,b]$, then:
\begin{equation*}
D^\alpha_t(f_1f_2)(t)=\sum\limits^\infty_{k=0}\left(\begin{array}{c}\alpha\\k\end{array}\right)D^{\alpha-k}_tf_1(t)\left
(\frac{d}{dt}\right)^kf_2(t),
\end{equation*} where $\left
(\ds\frac{d}{dt}\right)^k=\ds\frac{d}{dt}\circ...\circ\ds\frac{d}{dt};$

5. \begin{equation*}
\int^b_af_1(t)D^\alpha_tf_2(t)dt=-\int^b_af_2(t)^*D^\alpha_tf_1(t)dt;
\end{equation*}

6. If $f:[a,b]\rightarrow\mathbb R$ is analytical and $0\in(a,b)$
then

\begin{equation*}
f(t)=\sum\limits^\infty_{h=0}E_\alpha(t^h)D^{\alpha
h}_tf(t)|_{t=0},
\end{equation*} where $E_\alpha$ is the Mittag-Leffler function,
$E_\alpha(t^h)=\sum\limits^\infty_{h=0}\ds\frac{t^{\alpha
h}}{\Gamma(1+\alpha h)}.$

\end{proposition}

\subsection{The fractional osculator bundle}

\hspace{0.5cm} Let $\alpha\in(0,1)$ and M a n-dimensional
differentiable manifold. The parameterized curves on M,
$c_1,c_2:I\rightarrow M$, with $0\in I$, $c_1(0)=c_2(0)\in M$ have
a fractional contact $\alpha$ in $x_0$ if the relation
\begin{equation}
D^\alpha_t(f\circ c_1)|_{t=0}=D^\alpha_t(f\circ c_2)|_{t=o}
\end{equation} holds, for all $f\in {\cal F}(U)$ and $x_0\in U$, where $U$ is a local chart on
M.

Preceding equality (1) defines a relation of equivalence. The
classes $[c]^\alpha_{x_0}$ are called the fractional osculator
space in $x_0$, which will be denoted by
$Osc^{(\alpha)}_{x_0}(M)=T^{\alpha}_{x_0}(M)$.

Let $T^\alpha(M)=\bigcup\limits_{x_0\in M}T^\alpha_{x_0}(M)$ and
$\pi^\alpha:T^\alpha(M)\rightarrow M$, given by
$\pi^\alpha[c]^\alpha_{x_0}=x_0$. There is a differential
structure on $T^\alpha(M)$ and $(T^\alpha(M),\pi^\alpha,M)$ is a
differentiable bundle space.

If U is a local chart on M with $x_0\in U$ and $c:I\rightarrow M$
is a curve given by $x^i=x^i(t)$, $i=1..n$, $t\in I$, then
$[c]^\alpha_{x_0}$ is characterized by:
\begin{equation*}
x^i(t)=x^i(0)+\ds\frac{t^\alpha}{\Gamma(1+\alpha)}D^\alpha_tx^i|_{t=0},
\quad t\in(-\varepsilon,\varepsilon).
\end{equation*}

The coordinates of $[c]^\alpha_{x_0}$ on
$(\pi^\alpha)^{-1}(U)\subset T^\alpha_{x_0}(M)$ are
$(x^i,y^{i(\alpha)})$, where
\begin{equation*}
x^i=x^i(0),
y^{i(\alpha)}=\ds\frac{1}{\Gamma(1+\alpha)}D^\alpha_tx^i(t)|_{t=0},
\quad i=1..n.
\end{equation*}

From Proposition 2.1 and and the definition of $T^\alpha(M)$ we
have:

\begin{proposition}{\label{prop2}} 1. If $0<\alpha<\beta$ then
$T^{(\alpha)}(M)\subset T^{(\beta)}(M)$;

2. If $\underset{n\rightarrow\infty}{lim}\alpha_n=1$ then
$\underset{n\rightarrow\infty}{lim}T^{(\alpha_n)}(M)=T(M)$.
\end{proposition}

Let $(x^i)$, $i=1..n$ be the coordinate functions on U and
$(dx^i)_{i=1..n}$ be the base of 1-forms ${\cal{D}}^1(U)$ and
$\left(\ds\frac{\partial}{\partial x^i}\right )_{i=1..n}$ the base
of the vector fields ${\cal X}(U)$. For $f:U\rightarrow \mathbb R$
and $\alpha\in(0,1)$, the fractional derivative with respect to
$x^i$ is defined by:

\begin{equation}
\begin{split}
&D^\alpha_{x^i}f(x)=\\
&\!=\!\ds\frac{1}{\Gamma(1\!-\!\alpha)}\ds\frac{\partial}{\partial
x^i}\int^{x^i}_{a^i}\!\!\frac{f(x^1,\!...,\!x^{i\!-\!1},\!s,\!x^{i\!+\!1},\!...,\!x^n)\!-\!f(x^1,\!...,\!x^{i\!-\!1},a^i,x^{i\!+\!1},\!...,\!x^n)}{(x^i\!-\!s)^\alpha}ds
\end{split}
\end{equation} where $U_{ab}=\{x\in U, a^i\leq x^i\leq b^i,
i=1..n\}\subset U$.

From (2) we have:

\begin{proposition}{\label{prop3}} 1. If $f^i_1=(x^i)^\gamma$ then
$(D^\alpha_{x^i}f^i_1)(x)=\ds\frac{(x^i)^{j-\alpha}\Gamma(1+\alpha)}{\Gamma(1+\gamma-\alpha)};$

2. If $f^j_2=\ds\frac{(x^j)^\alpha}{\Gamma(1+\alpha)}$ then
$(D^\alpha_{x^i}f^j)(x)=\delta^j_i;$

3.
$D^\alpha_{x^i}(D^\alpha_{x^j}f)(x)=D^\alpha_{x^j}(D^\alpha_{x^i}f)(x)$,
$i,j=1..n.$
\end{proposition}

We consider the functions $(x^i)^\alpha\in{\cal{F}}(U)$ and
$d(x^i)^\alpha=\alpha(x^i)^{\alpha-1}dx^i\in{\cal{D}}^1(U)$,
$i=1..n$. The fractional exterior derivative is the operator
$d^\alpha:{\cal F}(U)\rightarrow{\cal D}^1(U)$ given by [1]:
\begin{equation*}
d^\alpha f=d(x^i)^\alpha D^\alpha_{x^i}(f).
\end{equation*}

Let $D^\alpha_{x^i}:{\cal{D}}^1(U)\rightarrow {\cal{D}}^1(U)$ be
the operator given by:
\begin{equation}
D^\alpha_{x^i}(a_jdx^j)=\sum\limits_{k=0}^\infty
\left(\begin{array}{c}\alpha\\k\end{array}\right)D^{\alpha-k}_{x^i}(a_j)\left
(\ds\frac{\partial}{\partial x^i}\right )^k(dx^j).
\end{equation}
From (3) and \begin{equation*}\left (\ds\frac{\partial}{\partial
x^i}\right )^k(dx^j)=d(\left (\ds\frac{\partial}{\partial
x^i}\right )^k(x^j))=0, k\geq 1
\end{equation*}
we obtain:
\begin{equation}
D^\alpha_{x^i}(a_jdx^j)=D_{x^i}^\alpha(a_j)dx^j.
\end{equation}

Let $d^\alpha:{\cal{D}}^1(U)\rightarrow {\cal{D}}^2(U)$ be the
operator given by:
\begin{equation}
d^\alpha(a_jdx^j)=d(x^i)^\alpha\wedge D^\alpha_{x^i}(a_jdx^j).
\end{equation}

From (4) and (5) we can deduce:
\begin{equation}
\begin{split}
&d^\alpha(a_jdx^j)=D^\alpha_{x^i}(a_j)d(x^i)^\alpha\wedge dx^j\\
&d^\alpha(b_jd(x^j)^\alpha)=D^\alpha_{x^i}(b_j)d(x^i)^\alpha\wedge
d(x^j)^\alpha.
\end{split}
\end{equation}

\begin{proposition}{\label{prop4}} Let $U$, $\ov U$, $U\cap\ov
U\neq\varnothing$ be two charts on M, $x\in U\cap \ov U$ and the
change of local chart given by:
\begin{equation}
{\ov x}^i={\ov x}^i(x^1,...,x^n),\ rang\left (\ds\frac{\partial
{\ov x}^i}{\partial x^j}\right )=n.
\end{equation}
With respect to (7) the following relations:
\begin{equation*}
\begin{split}
&d(x^i)^\alpha={\mathop{J}\limits^\alpha}{}^i_j(x,\ov x)d(\ov x^j)^\alpha\\
&D_{x^i}^\alpha={\mathop{J}\limits^\alpha}{}^j_i(\ov
x,x)D^\alpha_{\ov x^j}\\
&{\mathop{J}\limits^\alpha}{}^i_j(x,\ov
x)\cdot{\mathop{J}\limits^\alpha}{}^j_k(\ov x,x)=\delta^i_k\\
&{\mathop{J}\limits^\alpha}{}^i_j(x,\ov
x)=(x^i)^{\alpha-1}\ds\frac{\partial {x}^i}{\partial \ov
x^j}\ds\frac{1}{(\ov x^j)^{1-\alpha}},
\end{split}
\end{equation*} hold,
where
\begin{equation*}
{\mathop{J}\limits^\alpha}{}^i_j(x,\ov
x)=\ds\frac{1}{\Gamma(1+\alpha)}D^\alpha_{\ov x^j}(x^i)^\alpha.
\end{equation*}
\end{proposition}

Let ${\cal X}^\alpha(U)$ be the module of the fractional vector
fields generated by the operators $\{D^\alpha_{x^i}\}_{i=1..n}$. A
fractional field of vectors ${\mathop{X}\limits^{\alpha}}\in{\cal
X}^\alpha(U)$ has the form
${\mathop{X}\limits^\alpha}={\mathop{X}\limits^\alpha}{}^iD^\alpha_{x^i}$,
where ${\mathop{X^i}\limits^{\alpha}}\in{\cal F}(U), i=1..n$.
Under a change of local chart it changes by ${\mathop{\ov
X}\limits^\alpha}{}^i={\mathop{J}\limits^\alpha}{}^i_j(x,\ov
x)X^j$.

The fractional differential equation associated to the fractional
field of vectors ${\mathop{X}\limits^{\alpha}}$ is:
\begin{equation}
D^\alpha_tx^i(t)={\mathop{X}\limits^{\alpha}}{}^i(x(t)),\, i=1..n.
\end{equation}
The fractional differential equation (8) with initial conditions
has solutions [2].

\section{The fractional osculator bundle of higher order. Geometrical structures.}
\subsection{The fractional k-osculator bundle, $k\geq1$.}
\hspace{0.5cm} The parameterized curves on M,
$c_1,c_2:I\rightarrow M$, with $0\in I$, $c_1(0)=c_2(0)=x_0\in M$
have a fractional contact of k order in $x_0$ if for any $f\in
{\cal F}(U)$, the following relations:
\begin{equation*}
D^{\alpha a}_t(f\circ c_1)|_{t=0}=D^{\alpha a}_t(f\circ
c_2)|_{t=0}, \, a=1..k
\end{equation*} hold, where $x_0\in U$ and U is a local chart on M.

The classes $([c]^{\alpha a}_{x_0})_{a=1..k}$ are called the
fractional osculator space of k order and they will be denoted by
$Osc^{\alpha k}_{x_0}(M)={\mathop{E}\limits^{\alpha k}}{}_{x_0}$.

We consider ${\mathop{E}\limits^{\alpha k}}=\bigcup\limits_{x_0\in
M}{\mathop{E}\limits^{\alpha k}}{}_{x_0}$ and $\pi^{\alpha
k}_0:{\mathop{E}\limits^{\alpha k}}\rightarrow M$ given by
$\pi^{\alpha k}_0([c]^{\alpha a}_{x_0})_{a=1..k}=x_0$. There is a
differentiable structure on ${\mathop{E}\limits^{\alpha k}}$ and
$({\mathop{E}\limits^{\alpha k}}, \pi^{\alpha k}_0, M)$ is a
differentiable bundle.

If U is a local chart on M with $x_0\in U$ and $c:I\rightarrow M$
is a curve given by $x^i=x^i(t)$, $i=1..n$, $t\in I$, then a class
$([c]^{\alpha a}_{x_0})_{a=1..k}$ is given by the curve:
\begin{equation*}
x^i(t)=x^i(0)+\sum\limits^k_{a=1}\ds\frac{t^{\alpha
a}}{\Gamma(1+\alpha k)}D^{\alpha a}_tx^i(t)|_{t=0},\,
t\in(-\varepsilon, \varepsilon).
\end{equation*}

In $(\pi^{\alpha k}_0)^{-1}(U)\subset {\mathop{E}\limits^{\alpha
k}}$, the coordinates of $([c]^{\alpha a}_{x_0})_{a=1..k}$ are
$(x^i$,$y^{i(\alpha a)})$, $i=1..n$, $a=1..k$, where
\begin{equation*}
x^i=x^i(0),\, y^{i(\alpha a)}=\ds\frac{1}{\Gamma(1+\alpha
k)}D^{\alpha a}_tx^i(t),\, i=1..n,a=1..k.
\end{equation*}

\begin{proposition}{\label{prop5}} Let $U$, $\ov U$, $U\bigcap\ov
U\neq\emptyset$ be two charts on M and
\begin{equation}
\ov x^i=\ov x^i(x^1,...,x^n),\, i=1..n, det\left (
\ds\frac{\partial \ov x^i}{\partial x^j}\right)\neq0
\end{equation}
the coordinates transformation. The coordinates transformation on
$(\pi^{\alpha k}_0)^{-1}$ $(U\bigcap\ov U)$ are given by:
\begin{equation}
\begin{split}
&\ov y^{i(\alpha)}={\mathop{J}\limits^\alpha}{}^i_j(x,\ov
x)y^{j(\alpha)}\\
&\ds\frac{\Gamma(\alpha(a-1))}{\Gamma(\alpha)}\ov y^{i(\alpha
a)}=\Gamma(1+\alpha){\mathop{J}\limits^\alpha}{}^i_j({\ov
y}^{\alpha(a-1)},
x)y^{j(\alpha)}+\ds\frac{\Gamma(2\alpha)}{\Gamma(\alpha)}{\mathop{J}\limits^\alpha}{}^i_j(
\ov y^{\alpha(a-1)}, y^\alpha)y^{j(2\alpha)}+\\
&+...+\ds\frac{\Gamma(2\alpha)}{\Gamma(\alpha)}{\mathop{J}\limits^\alpha}{}^i_j(
\ov y^{\alpha(a-1)}, y^{\alpha b})y^{j(\alpha
b)}+...+\ds\frac{\Gamma(\alpha(a-1))}{\Gamma(\alpha)}y^{i(\alpha
a)},\\
&a=2..k,\, b=2..k,\, b\leq a,
\end{split}
\end{equation} where
\begin{equation*}
\begin{split}
&{\mathop{J}\limits^\alpha}{}^i_j(x,\ov
x)=D^\alpha_{\ov x^j}(x^i)\\
&{\mathop{J}\limits^\alpha}{}^i_j({\ov y}^{\alpha(a-1)},
y^{\alpha(b-1)})=D^\alpha_{y^{j(\alpha(b-1))}}
\ov y^{i(\alpha(a-1))}, \, a,b=2...k, b\leq a,\\
&{\mathop{J}\limits^\alpha}{}^i_j({\ov y}^{\alpha(a-1)},
x)=D^\alpha_{x^j} \ov y^{i(\alpha(a-1))}, \, i,j=1...n.
\end{split}
\end{equation*}
\end{proposition}

From the definition of fractional osculator bundle we can deduce
that if $\underset{n\rightarrow\infty}{lim}\alpha_n=1$ then
$\underset{n\rightarrow\infty}{lim}{\mathop{E}\limits^{\alpha_
n}}=E=Osc^k(M)$. The bundle space $(E, \pi, M)$ was defined and
studied in [6].

\subsection{Geometrical structures on ${\mathop{E}\limits^{\alpha k}}$.}

\hspace{0.5cm} Let $\pi^{\alpha k}_{\alpha
h}:{\mathop{E}\limits^{\alpha
k}}\rightarrow{\mathop{E}\limits^{\alpha h}}$, $h<k$, be the
projections given by:
\begin{equation*}
\pi^{\alpha k}_{\alpha h}(x,y^{(\alpha)},...,y^{(\alpha
k)})=(x,y^{(\alpha)},...,y^{(\alpha h)})
\end{equation*}
and the operator $d^\alpha\pi^{\alpha k}_{\alpha h}:{\cal
X}^\alpha({\mathop{E}\limits^{\alpha k}})\rightarrow{\cal
X}^\alpha({\mathop{E}\limits^{\alpha h}})$:
\begin{equation*}
d^\alpha\pi^{\alpha k}_{\alpha h}=\Gamma(1+\alpha)(d(x^i)^\alpha
D^\alpha_{x^i}+\sum\limits^{h}_{a=1}d(y^{i(\alpha a)})^\alpha
D^\alpha_{y^{i(\alpha a)}}), h<k,
\end{equation*}
where ${\cal X}^\alpha({\mathop{E}\limits^{\alpha k}})$ is the
module of fractional vector fields on ${\mathop{E}\limits^{\alpha
k}}$.

We consider ${\cal V}^{\alpha k}_{\alpha h}=Ker
d^\alpha\pi^{\alpha k}_{\alpha h}$, $h=0,1,...,k-1$ and its base
given by
$\{D^\alpha_{y^{i(\alpha(k+1))}},...,D^\alpha_{y^{i(\alpha
k)}}\}$, $i=1..n$. From the definition of ${\cal V}^{\alpha
k}_{\alpha h}$ we get:

\begin{equation*}
\begin{split}
{\cal V}^{\alpha k}_{\alpha (k-1)}\subset {\cal V}^{\alpha
k}_{\alpha (k-2)}\subset ...\subset {\cal V}^{\alpha
k}_{\alpha}\subset {\cal V}^{\alpha k}_{0}\\
d^\alpha\pi^{\alpha k}_{\alpha h}(D^\alpha _{x^i})=D^\alpha
_{x^i},\,  d^\alpha\pi^{\alpha k}_{\alpha h}(D^\alpha
_{y^{i(b)}})=D^\alpha _{y^{i(b)}}, \, b=1..h.
\end{split}
\end{equation*}
From Proposition 3.1 we obtain:
\begin{proposition}{\label{prop6}} Under the change of local chart (9), the operators $D^\alpha_{x^i}$,
$D^\alpha_{y^{i(\alpha a)}}$, $i=1..n$, $a=1..n$, change by:
\begin{equation}
\begin{split}
&D^\alpha_{x^i}=\sum\limits^k_{a=1}{\mathop{J}\limits^\alpha}{}^j_i(\ov
y^{(\alpha a)}, x)D^\alpha_{\ov y^{j(\alpha
a)}}+{\mathop{J}\limits^\alpha}{}^j_i(\ov x, x)D^\alpha_{\ov
x^{j}}\\
&D^\alpha_{y^{i(\alpha
a)}}=\sum\limits^k_{b=1}{\mathop{J}\limits^\alpha}{}^j_i(\ov
y^{\alpha b}, y^{\alpha a})D^\alpha_{\ov y^{j(\alpha b)}}, \quad
a=1..k.
\end{split}
\end{equation}
\end{proposition}
From Proposition 3.2 we can deduce that ${\cal V}^{\alpha
k}_{\alpha h}$ has geometrical character.

The following fractional fields of vectors:
\begin{equation}
\begin{split}
&{\mathop{\Gamma}\limits^\alpha}=y^{i(\alpha)}D^\alpha_{y^{i(\alpha
k)}}\\
&{\mathop{\Gamma}\limits^{2\alpha}}=\Gamma(1+\alpha)y^{i(\alpha)}D^\alpha_{y^{i(\alpha
(k-1))}}+\ds\frac{\Gamma(2\alpha)}{\Gamma(\alpha)}y^{i(2\alpha)}D^\alpha_{y^{i(\alpha
k)}}\\
&\dots\\
&{\mathop{\Gamma}\limits^{\alpha
k}}\!=\!\Gamma\!(\!1\!+\!\alpha)y^{i(\alpha)}D^\alpha_{y^{i(\alpha
)}}\!+\!\ds\frac{\Gamma(2\alpha)}{\Gamma(\alpha)}y^{i(2\alpha)}D^\alpha_{y^{i(2\alpha)}}\!+\!
\dots\!+\!\ds\frac{\Gamma(\alpha
(\!k\!-\!1\!))}{\Gamma(\alpha)}y^{i(\alpha
k)}D^\alpha_{y^{i(\alpha k)}}
\end{split}
\end{equation} are called Liouville fractional fields of vectors.

From (10) and (11) the fields ${\mathop{\Gamma}\limits^{\alpha
a}}$, $a=1..k$ have geometrical character and
${\mathop{\Gamma}\limits^{\alpha}}\in {\cal V}^{\alpha k}_{0}$,
${\mathop{\Gamma}\limits^{2\alpha}}\in {\cal V}^{\alpha
k}_{\alpha}$, ..., ${\mathop{\Gamma}\limits^{\alpha k}}\in {\cal
V}^{\alpha k}_{\alpha (k-1)}.$

The operators ${\mathop{J}\limits^{\alpha k}}:{\cal
X}^\alpha({\mathop{E}\limits^{\alpha k}})\rightarrow{\cal
X}^\alpha({\mathop{E}\limits^{\alpha k}})$ with the properties:
\begin{equation}
{\mathop{J}\limits^{\alpha
k}}(D^\alpha_{x^i})=D^\alpha_{y^{i(\alpha)}},\,
{\mathop{J}\limits^{\alpha k}}(D^\alpha_{y^{i(\alpha
a)}})=D^\alpha_{y^{i(\alpha (a+1))}},\, a=1..k-1,\,
{\mathop{J}\limits^{\alpha k}}(D^\alpha_{y^{i(\alpha k)}})=0
\end{equation}
is called $\alpha k$ fractional tangent structure.

From (11), (12) and (13) we have:

\begin{proposition}{\label{prop7}} $\alpha k$-fractional tangent
structure has the following properties:

1. ${\mathop{J}\limits^{\alpha k}}$ has a geometrical character;

2. rang$({\mathop{J}\limits^{\alpha k}})=kn$, \,
${\mathop{J}\limits^{\alpha k}}\circ ...\circ
{\mathop{J}\limits^{\alpha k}}=0;$

3. ${\mathop{J}\limits^{\alpha k}}({\mathop{\Gamma}\limits^{\alpha
k}})={\mathop{\Gamma}\limits^{\alpha (k-1)}}$, ...,
${\mathop{J}\limits^{\alpha
k}}({\mathop{\Gamma}\limits^{2\alpha}})={\mathop{\Gamma}\limits^{\alpha
}}$, ${\mathop{J}\limits^{\alpha
k}}({\mathop{\Gamma}\limits^{\alpha }})=0.$
\end{proposition}

The fractional field of vectors ${\mathop{S}\limits^{\alpha k}}\in
{\cal X}^\alpha({\mathop{E}\limits^{\alpha k}})$ is called $\alpha
k$-fractional spray if ${\mathop{J}\limits^{\alpha
k}}({\mathop{S}\limits^{\alpha
k}})={\mathop{\Gamma}\limits^{\alpha k}}$. From (12) and (13) we
obtain the form of ${\mathop{S}\limits^{\alpha k}}$:

\begin{equation}
\begin{split}
&{\mathop{S}\limits^{\alpha k}}=\Gamma(1+\alpha )y^{i(\alpha
)}D^{\alpha}_{x^i}
+\ds\frac{\Gamma(2\alpha )}{\Gamma(\alpha)}y^{i(2\alpha )}D^{\alpha}_{y^{i(\alpha)}}+...+\\
&+\ds\frac{\Gamma(\alpha (k-1))}{\Gamma(\alpha)}y^{i(\alpha
k)}D^{\alpha}_{y^{i(\alpha (k-1))}}-\ds\frac{\Gamma(\alpha
k)}{\Gamma(\alpha)}G^i(x,y^{(\alpha)},...,y^{(\alpha
k)})D^{\alpha}_{y^{i(\alpha k)}}.
\end{split}
\end{equation}

\begin{proposition}{\label{prop8}} The $\alpha k$-fractional spray uniquely defines the fractional differential equation given by:
\begin{equation*}
\ds\frac{1}{\Gamma(1+\alpha
k)}D^{\alpha(k+1)}_tx^i(t)+G^i(x,D^{\alpha}_tx,...,\Gamma(1+\alpha(k-1))D^{\alpha
k}_tx)=0.
\end{equation*}
\end{proposition}

Let ${\mathop{\cal N}\limits^{\alpha k}}$ be the submodule of
${\cal X}^\alpha({\mathop{E}\limits^{\alpha k}})$ so that:
\begin{equation*}
{\cal X}^\alpha({\mathop{E}\limits^{\alpha k}})|_{(\pi^{\alpha
k}_0)^{-1}(U)}={\cal N}^{\alpha k}\oplus{\cal V}^{\alpha
k}_0|_{(\pi^{\alpha k}_0)^{-1}(U)}.
\end{equation*}

The submodule ${\mathop{\cal N}\limits^{\alpha k}}$ is called
fractional nonlinear connection.

We consider $l^{\alpha k}:{\cal X}^{(\alpha)}(U)\rightarrow{\cal
X}^{\alpha}({\mathop{E}\limits^{\alpha k}})|_{(\pi^{\alpha
k}_0)^{-1}(U)}$, $l^{\alpha k}(D^\alpha_{x^i})=\Delta^{\alpha
k}_{x^i}$, $i=1..n$, where
\begin{equation}
\Delta^{\alpha
k}_{x^i}=D^\alpha_{x^i}-\sum\limits^k_{a=1}{\mathop{
N}\limits^{\alpha a}}{}^j_iD_{y^{i(\alpha a)}}, \, i=1..n.
\end{equation}
From (11) and (15) we obtain:
\begin{equation*}
\Delta^{\alpha k}_{{\ov
x}^i}={\mathop{J}\limits^{\alpha}}{}^j_i(x,\ov x)\Delta^{\alpha
k}_{x^i}.
\end{equation*}

The functions $({\mathop{N}\limits^{\alpha a}}{}^j_i)_{i,j=1..n;
a=1..k}$ are called the coefficients of fractional nonlinear
connection.

Let ${\mathop{\cal N_a}\limits^{\alpha k}}$ be the vertical
submodule given by:
\begin{equation}
{\mathop{{\cal N}}\limits^{\alpha k}}{}_0={\mathop{{\cal
N}}\limits^{\alpha k}},\, {\mathop{{\cal N}}\limits^{\alpha
k}}_a={\mathop{J}\limits^{\alpha k}}({\cal N}^{\alpha k}_{a-1}),\,
a=1..k-1,
\end{equation} where ${\mathop{{\cal N}}\limits^{\alpha k}}$ is
the submodule defined by fractional nonlinear connection.

From (16) we have:
\begin{equation}
{\cal X}^\alpha({\mathop{E}\limits^{\alpha k}})|_{(\pi^{\alpha
k}_0)^{-1}(U)}={\mathop{{\cal N}}\limits^{\alpha
k}}{}_0\oplus...\oplus{\mathop{{\cal N}}\limits^{\alpha
k}}{}_{k-1}\oplus{\cal V}^{\alpha k}_0|_{(\pi^{\alpha
k}_0)^{-1}(U)}, \, U\subset M.
\end{equation}

In what follows we will use a base adapted to the decomposition
(17). Then:
\begin{equation*}
\Delta^{\alpha k}_{y^{i(k)}}={\mathop{J}\limits^{\alpha
k}}(\Delta^{\alpha k}_{x^i}), \, \Delta^{\alpha k}_{y^{i(\alpha
a)}}={\mathop{J}\limits^{\alpha k}}(\Delta^{\alpha
k}_{y^{i(a-1)}}),\, a=2..k,\, D^{\alpha}_{y^{(\alpha)}},\, i=1..n,
\end{equation*}
with
\begin{equation*}
\Delta^{\alpha k}_{x^i}\in {\mathop{{\cal N}}\limits^{\alpha
k}}{}_{0},\, \Delta^{\alpha k}_{y^{i(\alpha a)}}\in {\mathop{{\cal
N}}\limits^{\alpha k}}{}_{a}, \, a=2..k-1, \,
D^{\alpha}_{y^{i(\alpha)}}\in {\cal V}^{\alpha k}_{0}, \, i=1..n
\end{equation*}
and
\begin{equation}
\begin{split}
&{\mathop{{\Delta}}\limits^{\alpha
k}}{}_{y^{i(\alpha)}}=D^{\alpha}_{y^{i(\alpha)}}-{\mathop{{
N}}\limits^{\alpha
}}{}_{i}^jD^{\alpha}_{y^{i(2\alpha)}}-...-{\mathop{{
N}}\limits^{\alpha(k-1)}}{}_{i}^jD^{\alpha}_{y^{j(\alpha k)}}\\
&{\mathop{{\Delta}}\limits^{\alpha
k}}{}_{y^{i(2\alpha)}}=D^{\alpha}_{y^{i(2\alpha)}}-...-{\mathop{{
N}}\limits^{\alpha(k-1)}}{}_{i}^jD^{\alpha}_{y^{j(\alpha k)}}\\
&\dots\\
&{\mathop{{\Delta}}\limits^{\alpha k}}{}_{y^{i(\alpha
(k-1))}}=D^{\alpha}_{y^{i(\alpha (k-1))}}-{\mathop{{
N}}\limits^{\alpha}}{}_{i}^jD^{\alpha}_{y^{j(\alpha k)}}.
\end{split}
\end{equation}
The dual base of (18) is:
\begin{equation}
\begin{split}
&{\mathop{{\delta}}\limits^{\alpha
k}}y^{i(\alpha)}=d(y^{i(\alpha)})^\alpha+{\mathop{{
M}}\limits^{\alpha
}}{}_{j}^id(x^j)^\alpha\\
&{\mathop{{\delta}}\limits^{\alpha
k}}y^{i(2\alpha)}=d(y^{i(2\alpha)})^\alpha+{\mathop{{
M}}\limits^{\alpha }}{}_{j}^id(y^{i(\alpha)})^\alpha+{\mathop{{
M}}\limits^{2\alpha
}}{}_{j}^id(x^{j})^\alpha\\
&\dots\\
&{\mathop{{\delta}}\limits^{\alpha k}}y^{i(\alpha
k)}=d(y^{i(\alpha k)})^\alpha+{\mathop{{ M}}\limits^{\alpha
}}{}_{j}^id(y^{i(\alpha(k-1))})^\alpha+\dots+{\mathop{{
M}}\limits^{\alpha k}}{}_{j}^id(x^{j})^\alpha,
\end{split}
\end{equation} where
\begin{equation*}
\begin{split}
&{\mathop{{ M}}\limits^{\alpha }}{}_{j}^i={\mathop{{
N}}\limits^{\alpha }}{}_{j}^i,\, {\mathop{{ M}}\limits^{2\alpha
}}{}_{j}^i={\mathop{{ N}}\limits^{2\alpha }}{}_{j}^i+{\mathop{{
N}}\limits^{\alpha }}{}_{k}^i{\mathop{{ N}}\limits^{\alpha
}}{}_{j}^h\\
&\dots\\
&{\mathop{{ M}}\limits^{\alpha k}}{}_{j}^i={\mathop{{
N}}\limits^{\alpha k}}{}_{j}^i+{\mathop{{ N_h}}\limits^{\alpha
(k-1) }}{\mathop{{ N}}\limits^{\alpha }}{}_{j}^h+\dots+{\mathop{{
N}}\limits^{\alpha}}{}_{h}^i{\mathop{{ N}}\limits^{\alpha (k-1)
}}{}_{j}^h.
\end{split}
\end{equation*} The functions ${\mathop{{ M}}\limits^{\alpha
a}}{}_{j}^i$, $a=1..k$ are called dual coefficients of fractional
nonlinear connection.

From (14) and (19) we obtain:
\begin{proposition}{\label{prop9}} A $\alpha k$-fractional spray ,
${\mathop{{S}}\limits^{\alpha k}}$, with the components $G^i(x$,
$y^{(\alpha)}$, $...$, $y^{(\alpha k)})$, determines a fractional
nonlinear connection with the dual coefficients given by:
\begin{equation*}
\begin{split}
&{\mathop{{ M}}\limits^{\alpha }}{}_{j}^i=D^\alpha_{y^{i(\alpha)}}G^i\\
&{\mathop{{ M}}\limits^{2\alpha
}}{}_{j}^i=\ds\frac{\Gamma(\alpha)}{\Gamma(2\alpha)}({\mathop{{
S}}\limits^{\alpha k}}({\mathop{{ M}}\limits^{\alpha
}}{}_{j}^i)+{\mathop{{M}}\limits^{\alpha }}{}^i_l{\mathop{{
M}}\limits^{\alpha }}{}^l_j)\\
&\dots\\
&{\mathop{{ M}}\limits^{\alpha
k}}{}_{j}^i=\ds\frac{\Gamma(\alpha(k-1))}{\Gamma(\alpha
k)}({\mathop{{ S}}\limits^{\alpha k}}({\mathop{{
M}}\limits^{\alpha(k-1)}}{}_{j}^i)+{\mathop{{M}}\limits^{\alpha
}}{}^i_l{\mathop{{ M}}\limits^{\alpha(k-1)}}{}^l_j).
\end{split}
\end{equation*}
\end{proposition}
We consider the adapted base given by (18) and the operator
${\mathop{\cal L}\limits^{\alpha k}}$ defined by:
\begin{equation}
\begin{split}
&{\mathop{{\cal L}}\limits^{\alpha k}}{}_{\Delta^{\alpha b}_{x^i}}({\mathop{\Delta}\limits^{\alpha k}}_{y^{j(\alpha a)}})=
{\mathop{{L}}\limits^{(\alpha k)}}{}_{ji}^h{\mathop{\Delta}\limits^{\alpha k}}_{y^{h(\alpha a)}}\\
&{\mathop{{\cal L}}\limits^{\alpha k}}{}_{\Delta^{\alpha
k}_{y^{i(\alpha b)}}}({\mathop{\Delta}\limits^{\alpha
k}}_{y^{j(\alpha a)}})=
{\mathop{{C}}\limits^{(\alpha b)}}{}_{ji}^h{\mathop{\Delta}\limits^{\alpha k}}_{y^{h(\alpha a)}},\, \alpha=0,1,...,k,\, b=1..k,\\
\end{split}
\end{equation} where $y^{i(0)}=x^i$. The coefficients $({\mathop{{L}}\limits^{(\alpha k)}}{}_{ji}^h, {\mathop{{C}}\limits^{(\alpha
b)}}{}_{ji}^h)$ are called the fractional coefficients of linear
connection N.

A distinguished fractional tensor field of type $(0,k)$ is given
by the following expression:
\begin{equation*}
{\mathop{{g}}\limits^{\alpha k}}={\mathop{{g}}\limits^{\alpha
k}}{}_{i_0i_1...i_k}d(y^{i_0(0)})^\alpha\otimes{\mathop{{\delta}}\limits^{\alpha
k}}y^{i_1(\alpha)}
\otimes\dots\otimes{\mathop{{\delta}}\limits^{\alpha
k}}y^{i_k(\alpha k)},
\end{equation*} where ${\mathop{{\delta}}\limits^{\alpha
k}}y^{i(\alpha a)}$, $a=0..k$ are given by (19) and
$y^{i(0)}=x^i.$

The covariant derivative with respect to fractional nonlinear
connection N of ${\mathop{{g}}\limits^{\alpha k}}$ is defined by:
\begin{equation*}
\begin{split}
&g^{\alpha k}_{i_0i_1...i_k|m}={\mathop{{\Delta}}\limits^{\alpha
k}}{}_{y^{m(0)}}(g^{\alpha
k}_{i_0i_1...i_k})-{\mathop{{L}}\limits^{(\alpha
k)}}{}_{i_0m}^j{\mathop{g}\limits^{\alpha k}}{}_{ji_1...i_k},\\
&g^{\alpha k}_{i_0i_1...i_k{\mathop{|}\limits^{(\alpha
b)}}m}={\mathop{{\Delta}}\limits^{\alpha k}}{}_{y^{m(\alpha
b)}}(g^{\alpha k}_{i_0i_1...i_k})-{\mathop{{C}}\limits^{\alpha
k}}{}_{i_1m}^hg^{\alpha
k}_{i_0h...i_k}-\dots-{\mathop{{C}}\limits^{\alpha
k}}{}_{i_km}^hg^{\alpha k}_{i_0...h}.
\end{split}
\end{equation*}

A fractional metric structure on ${\mathop{{E}}\limits^{\alpha
k}}$ is a fractional field of tensors of type $(0,2)$,
${\mathop{{g}}\limits^{\alpha k}}={\mathop{{g}}\limits^{\alpha
k}}{}_{ij}d(x^i)^\alpha\otimes d(x^j)^\alpha$, with
${\mathop{{g}}\limits^{\alpha
k}}{}_{ij}(x,y^{(\alpha)},...,y^{(\alpha k)})$, which is symmetric
and positively defined.

The fractional Sasaki lift of ${\mathop{{g}}\limits^{\alpha k}}$
is the fractional field of tensors given by:
\begin{equation*}
{\mathop{{G}}\limits^{\alpha k}}={\mathop{{g}}\limits^{\alpha
k}}_{ij}d(x^i)^\alpha\otimes
d(x^j)^\alpha+\sum\limits^k_{a=1}{\mathop{{g}}\limits^{\alpha
k}}_{ij}\delta y^{i(\alpha)}\otimes \delta y^{j(\alpha)}.
\end{equation*}

If:
\begin{equation*}
{\mathop{{g}}\limits^{\alpha k}}{}_{ij|m}=0,\,
{\mathop{{g}}\limits^{\alpha k}}{}_{ij{\mathop{|}\limits^{\alpha
k}}m}=0
\end{equation*} hold, then the fractional linear connection (20) is called
metrical.

\begin{proposition}{\label{prop8}} On ${\mathop{{E}}\limits^{\alpha
k}}$ there is a unique metrical fractional linear connection N
with respect to metrical structure ${\mathop{{G}}\limits^{\alpha
k}}$ with the property:
\begin{equation*}
{\mathop{{L}}\limits^{(\alpha
k)}}{}_{jl}^i={\mathop{{L}}\limits^{(\alpha k)}}{}_{lj}^i,\,
{\mathop{{C}}\limits^{(\alpha
a)}}{}_{jl}^i={\mathop{{C}}\limits^{(\alpha a)}}{}_{lj}^i,\,
a=1..k.
\end{equation*} The coefficients ${\mathop{{L}}\limits^{(\alpha
k)}}{}_{jl}^i$, ${\mathop{{C}}\limits^{(\alpha k)}}{}_{jl}^i$ have
the expressions:
\begin{equation*}
\begin{split}
{\mathop{{L}}\limits^{(\alpha
k)}}{}_{jl}^i=\frac{1}{2}{\mathop{{g}}\limits^{\alpha
k}}{}^{is}({\mathop{{\Delta}}\limits^{\alpha
k}}{}_{x^j}{\mathop{{g}}\limits^{\alpha
k}}{}_{sl}+{\mathop{{\Delta}}\limits^{\alpha
k}}{}_{x^l}{\mathop{{g}}\limits^{\alpha
k}}{}_{js}-{\mathop{{\Delta}}\limits^{\alpha
k}}{}_{x^s}{\mathop{{g}}\limits^{\alpha k}}{}_{jl})\\
{\mathop{{C}}\limits^{(\alpha
a)}}{}_{jl}^i=\frac{1}{2}{\mathop{{g}}\limits^{\alpha
k}}{}^{is}({\mathop{{\Delta}}\limits^{\alpha k}}{}_{y^{j(\alpha
a)}}{\mathop{{g}}\limits^{\alpha
k}}{}_{sl}+{\mathop{{\Delta}}\limits^{\alpha k}}{}_{y^{l(\alpha
a)}}{\mathop{{g}}\limits^{\alpha
k}}{}_{js}-{\mathop{{\Delta}}\limits^{\alpha k}}{}_{y^{s(\alpha
a)}}{\mathop{{g}}\limits^{\alpha k}}{}_{jl}),\, a=1..k.
\end{split}
\end{equation*}
\end{proposition}

\section{Lagrange space ${\mathop{{L}}\limits^{\alpha
k}}$. Applications.}

\subsection{The fractional Euler-Lagrange equation.}

\hspace{0.5cm} A fractional Lagrangian of k order, $k\in {\mathbb
N}^*$, on the differentiable manifold M, is a differentiable map
$L:{\mathop{{E}}\limits^{\alpha k}}\rightarrow{\mathbb R}$ on
${{\mathop{E}\limits^{\tilde{\alpha
k}}}}=\{(x,y^{(\alpha)},...,y^{(\alpha a)})\in
{\mathop{E}\limits^{{\alpha k}}}, rang||y^{i(\alpha)}||=1\}$.
Also, L is continuous in the points of
${\mathop{{E}}\limits^{\alpha k}}$ for which $y^{i(\alpha)}$ is
zero. Then,
\begin{equation*}
g_{ij}(x,y^{(\alpha)},...,y^{(\alpha
k)})=\frac{1}{2}D^\alpha_{y^{i(\alpha)}}D^\alpha_{y^{j(\alpha)}}L
\end{equation*} is d-fractional field of tensors on ${\mathop{{E}}\limits^{\alpha
k}}$. The Lagrangian L is regular if $rang(g_{ij})=n$ on
${\mathop{E}\limits^{\tilde {\alpha k}}}$.

Let $c:t\in [0,1]\rightarrow(x^i(t))\in M$ be a parameterized
curve so that $Im c\subset U$, where U is a chart on M. The
extension of curve c to ${\mathop{{E}}\limits^{\alpha k}}$,
${\mathop{{c}}\limits^{\alpha k}}$, is the following
differentiable map:
\begin{equation}
{\mathop{{c}}\limits^{\alpha
k}}:t\in[0,1]\rightarrow(x^i(t),y^{i(\alpha)}(t),...,y^{i(\alpha
k)}(t))\in {\mathop{{E}}\limits^{\tilde{\alpha k}}}.
\end{equation}
The action of L along the curve ${\mathop{{c}}\limits^{\alpha k}}$
is given by:
\begin{equation*}
I({\mathop{{c}}\limits^{\alpha
k}})=\int_0^1L(x(t),y^{(\alpha)}(t),...,y^{(\alpha k)}(t))dt.
\end{equation*}

Let $c_\varepsilon:t\in [0,1]\rightarrow(x^i(t,\varepsilon))\in M$
be the family of curves so that $Im c_\varepsilon\subset U$ and
$c_\varepsilon(0)=c(0)$, $y^{i(\alpha a)}(0)=y^{i(\alpha
a)}(1)=0$, $a=1..k-1$, where $\varepsilon$ is a sufficiently small
number in absolute value. The action on
${\mathop{{c_\varepsilon}}\limits^{\alpha k}}$ is:
\begin{equation*}
I({\mathop{{c_\varepsilon}}\limits^{\alpha
k}})=\int_0^1L(x_\varepsilon(t),y_\varepsilon^{(\alpha)}(t),...,y_\varepsilon^{(\alpha
k)}(t))dt.
\end{equation*} A necessary condition for $I({\mathop{{c}}\limits^{\alpha
k}})$ to be an extreme fractional value for
$I({\mathop{{c_\varepsilon}}\limits^{\alpha k}})$ is:
\begin{equation*}
D^\alpha_\varepsilon I({\mathop{{c_\varepsilon}}\limits^{\alpha
k}})|_\varepsilon=0.
\end{equation*}

By direct calculus we obtain:
\begin{proposition} The curve $c:t\in[0,1]\rightarrow(x^i(t))\in
M$ has the property that the action
$I({\mathop{{c}}\limits^{\alpha k}})$ is an extreme value of
$I({\mathop{{c_\varepsilon}}\limits^{\alpha k}})$ if $(x^i(t))$,
$i=1..n$ is a solution of fractional Euler-Lagrange equation:
\begin{equation}
D^\alpha_{x^i}L+\sum\limits^k_{a=1}(-1)^ad^{\alpha
a}_t(D^\alpha_{y^{i(\alpha a)}}L)=0,\quad i=1..n
\end{equation} where $d^{\alpha
a}_t=\sum\limits^a_{b=1}y^{i(\alpha b)}D^\alpha_{y^{i(\alpha
(b-1))}}$, $y^{i(0)}=x^i$, $a=1..h.$
\end{proposition}

A necessary condition for $I({\mathop{{c}}\limits^{\alpha k}})$ to
be an extreme value for
$I({\mathop{{c_\varepsilon}}\limits^{\alpha k}})$ is:
\begin{equation*}
\ds\frac{dI({\mathop{{c_\varepsilon}}\limits^{\alpha
k}})}{d\varepsilon}|_{\varepsilon=0}=0.
\end{equation*}

\begin{proposition}
The curve $c:t\in[0,1]\rightarrow(x^i(t))\in M$ has the property
that the action $I({\mathop{{c}}\limits^{\alpha k}})$ is an
extreme value of $I({\mathop{{c_\varepsilon}}\limits^{\alpha k}})$
if $(x^i(t))$, $i=1..n$ is a solution of fractional Euler-Lagrange
equation:
\begin{equation}
\ds\frac{\partial L}{\partial
x^i}+\sum\limits^k_{a=1}(-1)^ad^{\alpha }_t(\ds\frac{\partial
L}{\partial y^{i(\alpha a)}})=0,\, i=1..n,
\end{equation} where $d^{\alpha
}_t=\sum\limits^h_{a=1}y^{i(\alpha a)}D^\alpha_{y^{i(\alpha
(a-1))}}$, $y^{i(0)}=x^i$.
\end{proposition}

{\bf Example}. We consider the fractional differential equation:
\begin{equation}
\ds\frac{c\Gamma(1+\gamma)}{\Gamma(1+\gamma-\alpha)}x^{\gamma-\alpha}+
a_1\Gamma(1+2\alpha)y^{2\alpha}+a_2\Gamma(1+3\alpha)y^{3\alpha}+a_3\Gamma(1+4\alpha)y^{4\alpha}=0.
\end{equation} Equation (24) is the fractional Euler-Lagrange equation (22)
for the fractional Lagrange function:
\begin{equation*}
L=\ds\frac{c}{1+\gamma-\alpha}x^{\gamma}-
a_1\Gamma(1+2\alpha)(y^{\alpha})^\alpha+a_2\Gamma(1+3\alpha)(y^{2\alpha})^\alpha-a_3\Gamma(1+4\alpha)(y^{3\alpha})^\alpha.
\end{equation*} Equation (24) is the fractional Euler Lagrange
equation (23) for the fractional Lagrange function:
\begin{equation*}
\begin{split}
L&=\ds\frac{c\Gamma(1+\gamma)}{\Gamma(1+\gamma-\alpha)(\gamma-\alpha+1)}x^{\gamma-\alpha-1}-
\ds\frac{a_1}{2}\Gamma(1+2\alpha)(y^{\alpha})^2+\ds\frac{a_2}{2}\Gamma(1+3\alpha)(y^{2\alpha})^2-\\
&-\ds\frac{a_3}{2}\Gamma(1+4\alpha)(y^{3\alpha})^2.
\end{split}
\end{equation*}

Along the curve c we define the operators:
\begin{equation*}
\begin{split}
&{\mathop{{E}}\limits^{0}}_i=D^\alpha_{x^i}+\sum\limits^h_{a=1}(-1)^a\ds\frac{1}{\Gamma(1+\alpha
a)}d^{\alpha a}_t(D^\alpha_{y^{i(\alpha a)}}), \quad i=1..n\\
&{\mathop{{E}}\limits^{\alpha}}_i=\sum\limits^h_{a=1}(-1)^a\ds\frac{1}{\Gamma(1+\alpha
a)}d^{\alpha a}_t(D^\alpha_{y^{i(\alpha a)}}), \quad i=1..n\\
&\dots\\
&{\mathop{{E}}\limits^{\alpha
h}}_i=(-1)^k\ds\frac{1}{\Gamma(1+\alpha
k)}d^{\alpha k}_t(D^\alpha_{y^{i(\alpha k)}}), \quad i=1..n\\
\end{split}
\end{equation*} which have the property:
\begin{proposition} The operators ${\mathop{{E}}\limits^{0}}_i(L)$, ...,${\mathop{{E}}\limits^{\alpha
k}}_i(L)$, $i=1..n$, are d-fractional fields of covectors for any
differentiable Lagrangian of order $\alpha k$,
${\mathop{{L}}\limits^{\alpha k}}$, along the extension
${\mathop{{c}}\limits^{\alpha k}}$ of curve c.
\end{proposition} d-fields of covectors ${\mathop{{E}}\limits^{\alpha}}_i(L)$, ...,${\mathop{{E}}\limits^{\alpha
k}}_i(L)$ are called Craig and Synge d-fractional fields of
covectors.
\begin{proposition} 1. d-fractional field of covectors,
${\mathop{{E_i}}\limits^{\alpha(k-1)}}(L)$, has the form:
\begin{equation*}
{\mathop{E_i}\limits^{\alpha
(k-1)}}(L)=(-1)^{k-1}\ds\frac{1}{\Gamma(1+\alpha
(k-1))}(D^\alpha_{y^{i(\alpha
(k-1))}}L-{\mathop{\Gamma}\limits^{\alpha}}(D^\alpha_{y^{i(\alpha
(k))}}L)-g_{ij}y^{j(\alpha(k+1))}),
\end{equation*} $i=1..n,$ where ${\mathop{{\Gamma}}\limits^{\alpha}}$ is
given (12).

2. The system of fractional differential equations:
\begin{equation*}
g^{ij}{\mathop{{E}}\limits^{\alpha(k-1)}}{}_j(L)=0,\, i=1..n
\end{equation*} determines a $\alpha k$
fractional spray ${\mathop{{S}}\limits^{\alpha k}}$ on the curve
${\mathop{{c}}\limits^{\alpha k}}$, given by (21):
\begin{equation*}
{\mathop{G}\limits^{\alpha
ki}}{}=\ds\frac{\Gamma(\alpha)}{\Gamma(1+\alpha
k)\Gamma(1+\alpha)}g^{ij}[{\mathop{\Gamma}\limits^{\alpha}}(D^\alpha_{y^{j(\alpha
k)}}L)-D^\alpha_{y^{j(\alpha (k-1))}}], \, i,j=1..n.
\end{equation*}
\end{proposition}

\subsection{The prolongation of Riemann, Finsler and Lagrange fractional structures to fractional bundle of k order.}

\hspace{0.5cm} The pair ${\mathop{{\cal
R}}\limits^{\alpha}}=(M,{\mathop{g}\limits^{\alpha}})$ is called
Riemann fractional structure, where M is a differentiable manifold
of n dimension and
${\mathop{g}\limits^{\alpha}}=({\mathop{g}\limits^{\alpha}}{}_{ij})$
is a fractional field of tensors, which means that under a change
of local chart on M, the system of functions
${\mathop{g}\limits^{\alpha}}{}_{ij}$ change by:
\begin{equation*}
{\mathop{g}\limits^{\alpha}}_{ij}(\overline
x)={\mathop{J}\limits^{\alpha}}_{i}^l(\ov x,
x){\mathop{J}\limits^{\alpha}}_{j}^h(\ov x,
x){\mathop{g}\limits^{\alpha}}_{lh}(x)
\end{equation*} and
${\mathop{g}\limits^{\alpha}}_{ij}={\mathop{g}\limits^{\alpha}}_{ji}$
with $({\mathop{g}\limits^{\alpha}}_{ij})$ is positively defined.
The fractional Christofel symbols
${\mathop{\gamma^l}\limits^{\alpha}}_{ij}$ of
${\mathop{g}\limits^{\alpha}}$ are:
\begin{equation*}
{\mathop{\gamma}\limits^{\alpha}}{}_{ij}^l=\frac{1}{2}
{\mathop{g}\limits^{\alpha}}{}^{ls}(D^\alpha_{x^i}{\mathop{g}\limits^{\alpha}}_{sj}+
D^\alpha_{x^j}{\mathop{g}\limits^{\alpha}}_{is}-D^\alpha_{x^s}{\mathop{g}\limits^{\alpha}}_{ij}).
\end{equation*}
The prolongation of ${\mathop{g}\limits^{\alpha}}$ to
${\mathop{E}\limits^{\alpha k}}$ is the fractional field of
tensors ${\mathop{g}\limits^{\alpha k}}$ with the property:
\begin{equation*}
({\mathop{g}\limits^{\alpha k}}\circ \pi^{\alpha
k}_0)(x,y^{(\alpha)},..., y^{(\alpha
k)})={\mathop{g}\limits^{\alpha}}(x),\, \forall
(x,y^{(\alpha)},..., y^{(\alpha k)})\in (\pi^{\alpha
k}_0)^{-1}(U).
\end{equation*}

\begin{proposition}There are fractional nonlinear connections ${\mathop{N}\limits^{\alpha
k}}$ on ${\mathop{E}\limits^{\alpha k}}$ which are determined only
by ${\mathop{g}\limits^{\alpha}}$. One of them is:
\begin{equation}
\begin{split}
&{\mathop{M}\limits^{\alpha}}{}_j^i={\mathop{\gamma}\limits^{\alpha}}{}_{jm}^iy^{(\alpha)m},\\
&{\mathop{M}\limits^{2\alpha}}{}_j^i=\ds\frac{\Gamma(\alpha)}{\Gamma(2\alpha)}({\mathop{\Gamma}\limits^{\alpha}}({\mathop{M}\limits^{\alpha}}{}_j^i)+
{\mathop{M}\limits^{\alpha}}{}_h^i{\mathop{M}\limits^{\alpha}}{}_j^h),\\
&\dots\\
&{\mathop{M}\limits^{\alpha
k}}{}_j^i=\ds\frac{\Gamma(\alpha(k-1))}{\Gamma(\alpha
k)}({\mathop{\Gamma}\limits^{\alpha}}({\mathop{M}\limits^{\alpha(k-1)}}{}_j^i)+
{\mathop{M}\limits^{\alpha}}{}_h^i{\mathop{M}\limits^{\alpha(k-1)}}{}_j^h).
\end{split}
\end{equation}
\end{proposition} For $k=1$ the coefficients of fractional
nonlinear connection ${\mathop{N}\limits^{\alpha}}$ on
${\mathop{E}\limits^{\alpha}}$ are
${\mathop{M}\limits^{\alpha}}{}_j^i=\gamma^i_{jh}y^{h(\alpha)}$
and for $k=2$ the coefficients are:
\begin{equation*}
{\mathop{M}\limits^{\alpha}}{}_j^i=\gamma^i_{jh}y^{h(\alpha)},
{\mathop{M}\limits^{2\alpha}}{}_j^i=\gamma^i_{jh}y^{h(2\alpha)}+\frac{\Gamma(\alpha)}
{\Gamma(2\alpha)}(D^\alpha_{x^h}\gamma^i_{jp}+\gamma^i_{hl}\gamma^l_{jp})y^{h(\alpha)}y^{p(2\alpha)}.
\end{equation*}

The pair ${\mathop{{\cal F}}\limits^{\alpha}}=(M,
{\mathop{F}\limits^{\alpha}})$ is called Finsler fractional
structure, where M is a differentiable manifold of dimension n,
and
${\mathop{F}\limits^{\alpha}}:{\mathop{E}\limits^{\alpha}}\rightarrow
M$ is called fundamental fractional function. We consider the
prolongation of ${\mathop{F}\limits^{\alpha}}$,
${\mathop{F}\limits^{\alpha k}}:{\mathop{E}\limits^{\alpha
k}}\rightarrow {\mathbb R}$ given by:
\begin{equation*}
({\mathop{F}\limits^{\alpha k}}\circ \pi^{\alpha
k}_\alpha)(x,y^{(\alpha)},...,y^{(\alpha
k)})={\mathop{F}\limits^{\alpha}}(x,y^{(\alpha)})
\end{equation*} and the prolongation of fundamental fractional
tensor:
\begin{equation*}
{\mathop{\gamma}\limits^{\alpha
}}{}_{ij}(x,y^{(\alpha)})=\frac{1}{2}D^\alpha_{y^{i(\alpha)}}(D^\alpha_{y^{j(\alpha)}}{\mathop{F^2}\limits^{\alpha
}})
\end{equation*} given by:
\begin{equation*}
({\mathop{\gamma}\limits^{\alpha k}}{}_{ij}\circ \pi^{\alpha
k}_0)(x,y^{(\alpha)},...,y^{(\alpha
k)})={\mathop{\gamma}\limits^{\alpha }}{}_{ij}(x,y^{\alpha}).
\end{equation*} Let ${\mathop{\gamma}\limits^{\alpha }}{}_{ij}^h(x,y^{(\alpha)})$ be the Christoffel symbols of $({\mathop{\gamma}\limits^{\alpha
}}{}_{ij})$, given by:
\begin{equation*}
{\mathop{\gamma}\limits^{\alpha}}{}_{ij}^h(x,y^{(\alpha)})=\frac{1}{2}{\mathop{\gamma}\limits^{\alpha}}{}^{ls}
(D^\alpha_{x^i}{\mathop{\gamma}\limits^{\alpha}}{}_{ij}+D^\alpha_{x^j}{\mathop{\gamma}\limits^{\alpha}}{}_{is}-
D^\alpha_{x^s}{\mathop{\gamma}\limits^{\alpha}}{}_{ij}),
\end{equation*} where
$({\mathop{\gamma}\limits^{\alpha}}{}^{ls})=({\mathop{\gamma}\limits^{\alpha}}{}_{ls})^{-1}$.

The coefficients of nonlinear fractional connection (fractional
Cartan coefficients) are:
\begin{equation*}
{\mathop{G}\limits^{\alpha}}{}_{j}^i=\frac{1}{2}D^\alpha_{y^{j(\alpha)}}({\mathop{\gamma}\limits^{\alpha}}{}_{pm}^iy^{p(\alpha)}y^{m(\alpha)}).
\end{equation*}
\begin{proposition} There is a nonlinear fractional connection on
${\mathop{E}\limits^{\tilde{\alpha k}}}={\mathop{E}\limits^{\alpha
k}}\diagdown \{0\}=\{(x,y^{(a)},...,y^{(\alpha a)})\in
{\mathop{E}\limits^{\alpha k}}, rang||y^{i(\alpha)}||=1\}$ which
only depends on the fundamental fractional function
${\mathop{F}\limits^{\alpha}}$ of the fractional Finsler space.
One of them has the dual coefficients given by:
\begin{equation*}
\begin{split}
&{\mathop{M}\limits^{\alpha}}{}_j^i={\mathop{G}\limits^{\alpha}}{}_{j}^i,\\
&{\mathop{M}\limits^{2\alpha}}{}_j^i=\ds\frac{\Gamma(\alpha)}{\Gamma(2\alpha)}({\mathop{\Gamma}\limits^{\alpha}}({\mathop{M}\limits^{\alpha}}{}_j^i)+
{\mathop{G}\limits^{\alpha}}{}_m^i{\mathop{M}\limits^{\alpha}}{}_j^m),\\
&\dots\\
&{\mathop{M}\limits^{\alpha
k}}{}_j^i=\ds\frac{\Gamma(\alpha(k-1))}{\Gamma(\alpha
k)}({\mathop{\Gamma}\limits^{\alpha}}({\mathop{M}\limits^{\alpha(k-1)}}{}_j^i)+
{\mathop{G}\limits^{\alpha}}{}_m^i{\mathop{M}\limits^{\alpha(k-1)}}{}_j^m).
\end{split}
\end{equation*}
\end{proposition}

The pair ${\mathop{{\cal L}}\limits^{\alpha}}=(M,
{\mathop{L}\limits^{\alpha}})$ is called Lagrange fractional
structure, where
${\mathop{L}\limits^{\alpha}}:{\mathop{E}\limits^{\alpha}}\rightarrow
{\mathbb R}$.

The prolongation of ${\mathop{L}\limits^{\alpha}}$ to
${\mathop{E}\limits^{\alpha k}}$ is defined by:
\begin{equation*}
({\mathop{L}\limits^{\alpha k}}\circ \pi^{\alpha
k}_\alpha)(x,y^{(\alpha)},..., y^{(\alpha
k)})={\mathop{L}\limits^{\alpha}}(x,y^{(\alpha)})
\end{equation*} and the prolongation of fundamental tensor:
\begin{equation*}
{\mathop{g}\limits^{\alpha}}{}_{ij}(x,y^{(\alpha)})=\frac{1}{2}D^\alpha_{y^{i(\alpha)}}D^\alpha_{y^{j(\alpha)}}{\mathop{L}\limits^{\alpha
}}(x,y^{(\alpha)})
\end{equation*} to ${\mathop{E}\limits^{\alpha k}}$ is:
\begin{equation*}
({\mathop{g}\limits^{\alpha k}}_{ij}\circ \pi^{\alpha
k}_\alpha)(x,y^{(\alpha)},..., y^{(\alpha
k)})={\mathop{g}\limits^{\alpha}}_{ij}(x,y^{(\alpha)}).
\end{equation*}

Considering the integral action
$I({\mathop{c}\limits^{\alpha}})=\int_0^1L(x(t),y^{(\alpha)}(t))dt$
on a parameterized curve $c$, the fractional Euler-Lagrange
equations (22) are:
\begin{equation*}
y^{i(2\alpha)}={\mathop{G}\limits^{\alpha}}{}^i(x,y^{(\alpha)}),\,
i=1..n
\end{equation*} where
\begin{equation*}
{\mathop{G}\limits^{\alpha}}{}^i(x,y^{(\alpha)})=g^{im}(D^{\alpha}_{y^{m(\alpha)}}D^\alpha_{x^j}y^{j(\alpha)}-D^\alpha_{x^m}L),\,
i=1..n
\end{equation*} $(g^{im})=(g_{im})^{-1}$, $g_{im}=D^{\alpha}_{y^{i(\alpha)}}D^\alpha_{y^{m(\alpha)}}L.$

The functions: \begin{equation}
{\mathop{G}\limits^{\alpha}}{}^i_j(x,y^{(\alpha)})=D^{\alpha}_{y^{j(\alpha)}}G^i(x,y^{(\alpha)}),\,
i,j=1..n
\end{equation} are the first dual coefficients of fractional
nonlinear connection ${\mathop{{\cal N}}\limits^{\alpha}}$ on
${\mathop{E}\limits^{\tilde\alpha}}$ which only depends on the
fundamental function ${\mathop{L}\limits^{\alpha}}$ of the
Lagrange space ${\mathop{\cal L}\limits^{\alpha}}$.

We obtain the result: \begin{proposition} There are nonlinear
fractional connections on ${\mathop{E}\limits^{\tilde{\alpha k}}}$
which only depend on the fundamental function
${\mathop{L}\limits^{\alpha}}$ of the Lagrange space
${\mathop{{\cal L}}\limits^{\alpha}}$. One of them has the dual
coefficients (25), where ${\mathop{G}\limits^{\alpha}}{}_j^i$ are
given by (26).
\end{proposition}

\section{Conclusions.}

\hspace{0.5cm} The study conducted in this paper takes into
account the geometrical character of the introduced objects. In
the case $M={\mathbb R}$, using the methods from this paper, the
information concerning fractional differential systems which
describe concrete processes was obtained in [3]. The results from
the present paper will permit the study of other geometrical
objects which are described with the help of the fractional
derivative.

\end{document}